\documentclass{jams-l}
\usepackage{amsmath}
\usepackage{amsthm}
\usepackage{amsfonts}
\usepackage{amssymb}

\def\R{{\hbox{\bf R}}}

\def\Q{{\hbox{\bf Q}}}

\def\P{{\hbox{\bf P}}}
\def\E{{\hbox{\bf E}}}

\font \roman = cmr10 at 10 true pt

\def\be#1{ \begin{equation}\label{#1} }

\def\bas{\begin{align*}}
\def\eas{\end{align*}}
\def\bi{\begin{itemize}}
\def\ei{\end{itemize}}

\def\Gr{{\hbox{\roman Gr}}}

\def\Z{{\hbox{\bf Z}}}

\def\eps{\varepsilon}
\def \endprf{\hfill  {\vrule height6pt width6pt depth0pt}\medskip}
\def\emph#1{{\it #1}}
\def\textbf#1{{\bf #1}}

\def\ep{{\epsilon}}

\theoremstyle{plain}
  \newtheorem{theorem}[subsection]{Theorem}
  \newtheorem{conjecture}[subsection]{Conjecture}

  \newtheorem{proposition}[subsection]{Proposition}
  
  \newtheorem{lemma}[subsection]{Lemma}

\theoremstyle{remark}
  \newtheorem{remark}[subsection]{Remark}

\theoremstyle{definition}
  \newtheorem{definition}[subsection]{Definition}

\numberwithin{equation}{section}

\include{psfig}

\begin{document}

\title{On the singularity probability of random Bernoulli matrices}

\author{Terence Tao}
\address{Department of Mathematics, UCLA, Los Angeles CA 90095-1555}
\email{tao@math.ucla.edu}
\thanks{T. Tao is a Clay Prize Fellow and is supported by a grant from the Packard Foundation.}

\author{Van Vu}
\address{Department of Mathematics, Rutgers University, Piscataway, NJ 08854-8019}
\email{vanvu@ucsd.edu}
\thanks{V. Vu is an A. Sloan Fellow and is supported by an NSF Career Grant.}

\begin{abstract} Let $n$ be a large integer and $M_n$ be a random $n$ by
$n$ matrix whose entries are i.i.d. Bernoulli random variables
(each entry is $\pm 1$ with probability $1/2$). We show that the
probability that $M_n$ is singular is at most $(3/4 +o(1))^n$,
improving an earlier estimate of Kahn, Koml\'os and Szemer\'edi \cite{kks}, as well as earlier work by the authors \cite{TV1}.
The  key new ingredient is the applications of Freiman type
inverse theorems and other tools  from  additive combinatorics.
\end{abstract}


\maketitle

\section{Introduction} Let $n$ be a large integer, and $M_n$ be a random $n$ by
$n$ matrix whose entries are i.i.d. Bernoulli random variables
(each entry is $\pm 1$ with probability $1/2$). In this paper we
consider the well known question of estimating $P_n := \P(\det(M_n)=0)$, the
probability that $M_n$ is singular.

By considering the event that two columns or two rows of $M_n$ are
equal (up to sign), it is clear that

$$ P_n \geq (1 + o(1)) n^2 2^{1-n}.$$
\noindent
 Here and in the sequel we  use the asymptotic notation under the assumption that
  $n$ tends to infinity.
It has been conjectured by many researchers that in fact

\begin{conjecture} \label{conj1}
$$ P_n = (1 + o(1)) n^2 2^{1-n}= (\frac{1}{2}+ o(1))^n.$$
\end{conjecture}

In a breakthrough paper \cite{kks},  Kahn, Koml\'os and
Szemer\'edi proved that

\begin{equation} \label{KKS} P_n = O( .999^n ). \end{equation}

  In a recent paper \cite{TV1} we  improved
  this bound slightly to $O(.958^n)$. The main result of this paper is the following
  more noticeable improvement.

\begin{theorem}\label{main}
$$ P_n \le  (\frac{3}{4} + o(1))^n.$$
\end{theorem}

The proof of this theorem uses several ideas from \cite{kks, TV1}.
The key new ingredient is the applications of  Freiman type
inverse theorems and other tools  from  additive combinatorics.

In the next section, we present the  initial steps of the proof of
Theorem \ref{main}. These steps will reveal the origin of the
bound $(\frac{3}{4} + o(1))^n$ and also provide an approach to the
conjectured bound $(1+o(1))n^2 2^{1-n}$.

\section{Reduction to medium combinatorial dimension}

We begin the proof of Theorem \ref{main}.  We first make a
convenient (though somewhat artificial) finite field reduction.
Observe that any $n \times n$ matrix with entries $\pm 1$ can have
determinant at most $n^{n/2}$ (since the magnitude of the
determinant equals the volume of a $n$-dimensional parallelopiped,
all of whose edges have length $\sqrt{n}$); this is also known as
\emph{Hadamard's inequality}. Thus if we introduce a finite field
$F := \Z/p\Z$ of prime order
\begin{equation}\label{F-size}
|F| = p > n^{n/2},
\end{equation}
then $P_n$ is also the probability that a random $n
\times n$ matrix with entries $\pm 1$ (now thought of as elements
of $F$) is singular in the finite field $F$.  Henceforth we fix
the finite field $F$, and shall work over the finite field $F$
instead of over $\R$.  In particular, linear algebra terminology
such as ``dimension'', ``rank'', ``linearly independent'',
``subspace'', ``span'', etc. will now be with respect to the field
$F$.  One advantage of this discretized setting is that all
collections of objects under consideration will now automatically
be finite, and one can perform linear change of variables over $F$ without
having to worry about Jacobian factors.  Note that while one could use Bertrand's postulate
to fix $|F|$ to be comparable to $n^{n/2}$, we will not need to do so in this paper.  Indeed it will be more convenient to take $|F|$ to
be extremely large (e.g. larger than $\exp(\exp(Cn))$) to avoid any ``torsion'' issues.

Now let $\{-1,1\}^n \subset F^n$ be the discrete unit cube in
$F^n$.  We let $X$ be the random variable taking values in
$\{-1,1\}^n$ which is distributed uniformly on this cube (thus
each element of $\{-1,1\}^n$ is attained with probability
$2^{-n}$).  Let $X_1, \ldots, X_n \in \{-1,1\}$ be $n$ independent
samples of $X$.  Then
$$ P_n := \P( X_1, \ldots, X_n \hbox{ linearly dependent} ).$$
For each linear subspace $V$ of $F^n$, let $A_V$ denote the event that
$X_1,\ldots,X_n$ span $V$.
Let us call a space $V$ \emph{non-trivial} if it is spanned by the set $V \cap \{-1,1\}^n$.
Note that $\P(A_V) \neq 0$ if and only if $V$ is non-trivial.
Since every collection of $n$ linearly dependent
vectors in $F^n$ will span exactly one proper subspace $V$ of
$F^n$, we have
\begin{equation}\label{pn-spanned}
 P_n = \sum_{V \hbox{ a proper non-trivial subspace of } F^n} \P( A_V ).
\end{equation}

In \cite[Lemma 5.1]{TV1} (see also \cite{kks}), we showed that the dominant contribution to this sum came from the hyperplanes:
$$ P_n = 2^{o(n)} \sum_{V \hbox{ a non-trivial hyperplane in } F^n} \P( A_V ).$$
Thus it will suffice to show that
$$ \sum_{V \hbox{ a non-trivial hyperplane in } F^n} \P( A_V ) \leq (3/4 + o(1))^n.$$
We remark that the arguments below also extend to control the lower-dimensional spaces directly, but it will be somewhat simpler technically
to work just with hyperplanes.

As in \cite{kks}, the next step is to partition the non-trivial hyperplanes $V$ into a number of classes, depending on
a quantity which we shall call the \emph{combinatorial dimension}.

\begin{definition}[Combinatorial dimension]  Let $D := \{ d_\pm \in \Z/n:
1  \leq d_\pm \leq n\}$.
 For any $d_\pm \in D$, we define the \emph{combinatorial Grassmannian}
$\Gr(d_\pm)$ to be the set of all non-trivial
hyperplanes $V$ in $F^n$ with
\begin{equation}\label{discrete-dim}
2^{d_\pm - 1/n} < |V \cap \{-1,1\}^n| \leq 2^{d_\pm}.
\end{equation}
We will
refer to $d_\pm$ as the \emph{combinatorial dimension} of
$V$.
\end{definition}

\begin{remark} Following \cite{kks}, we have gradated the combinatorial dimension in steps of $1/n$ rather than steps of 1.  This is because
we will need to raise probabilities such as $\P(X \in V)$ to the power $n$, and thus errors of $2^{O(1/n)}$ in
will be acceptable while errors of $2^{O(1)}$ can be
more problematic.  This of course increases the number of possible
combinatorial dimensions from $O(n)$ to $O(n^2)$, but this will
have a negligible impact on our analysis.
\end{remark}

It thus suffices to show that
\begin{equation}\label{pn-spanned-2}
 \sum_{d_\pm \in D} \sum_{V \in \Gr(d_\pm)} \P( A_V ) \leq (\frac{3}{4} + o(1))^n.
\end{equation}
It is therefore of interest to understand the size of the combinatorial Grassmannians $\Gr(d_\pm)$
and of the probability of the events $A_V$ for hyperplanes $V$ in those Grassmannians.

There are two easy cases, one when $d_\pm$ is fairly small and one where $d_\pm$ is fairly large, both of which were
treated in \cite{kks}, and which we present below.

\begin{lemma}[Small combinatorial dimension estimate]\label{kks-small}\cite{kks} Let $0 < \alpha < 1$ be arbitrary.  Then
$$
\sum_{d_\pm \in D: 2^{d_\pm-n} \leq \alpha^n} \sum_{V \in \Gr(d_\pm)} \P( A_V )
\leq n\alpha^n.$$
\end{lemma}

\begin{proof} Observe that if $X_1, \dots, X_n$ span
$V$, then there are $n-1$ vectors among the $X_i$ which already
span $V$. By symmetry, we thus have

\begin{equation}\label{V-span}
 \P( A_V ) = \P( X_1,\ldots,X_n \hbox{ span } V)
 \leq n \P( X_1,\ldots,X_{n-1} \hbox{ span } V) \P( X \in V
 ).
 \end{equation}
On the other hand, if $V \in \Gr(d_\pm)$ and $2^{d_\pm-n} \leq \alpha^n$,
then $\P (X \in V) \le \alpha^n$ thanks to \eqref{discrete-dim}.  Thus we have
$$  \P( A_V )
 \leq n \alpha^n \P( X_1,\ldots,X_{n-1} \hbox{ span } V).$$
Since $X_1,\ldots,X_{n-1}$ can span at most one space $V$, the claim follows.
\end{proof}

\begin{lemma}[Large combinatorial dimension estimate]\label{kks-large}\cite{kks}  We have
$$
\sum_{d_\pm \in D: 2^{d_\pm-n} \geq 100 /\sqrt{n} } \,\,\,\,
\sum_{V \in \Gr(d_\pm)} \P( A_V ) \leq (1 + o(1))n^2
2^{-n}.$$
\end{lemma}

\begin{proof}(Sketch)  For this proof it is convenient to work back in Euclidean space $\R^n$ instead of in the finite field model $F^n$; observe that
one can identify non-trivial hyperplanes in $F^n$ with non-trivial hyperplanes in $\R^n$ without affecting the combinatorial dimension.
Let $V \subset \R^n$ be a hyperplane in $\Gr(d_\pm)$ with $2^{d_\pm-n} \geq \frac{100}{\sqrt{n}}$.  Then by \eqref{discrete-dim} we have
$$ \P( X \in V) \geq 50 / \sqrt{n}$$
(say). Since $V$ is spanned by elements $V \cap \{-1,1\}^n$ which lie in the integer lattice $\Z^n$, the orthogonal complement $V^\perp \subset \R^n$
contains at least
one integer vector $a = (a_1,\ldots,a_n) \in \Z^n \backslash \{0\}$ which is not identically zero.  In particular we see that
$$ \P( a_1 X_1 + \ldots + a_n X_n ) \geq 50 / \sqrt{n}.$$
But by Erd\"os's Littlewood-Offord inequality (see \cite{erdos}) this forces at most $n/2$ (say) of the coefficients $a_1,\ldots,a_n$
to be non-zero, if $n$ is sufficiently large.  Thus
\begin{align*}
&  \sum_{d_\pm \in D: 2^{d_\pm-n} \geq \frac{100}{\sqrt{n}}}
\,\,\,\,\, \sum_{V \in \Gr(d_\pm)} \P( A_V ) \\ &\leq \P(
a \cdot X_1 = \ldots = a \cdot X_n = 0\hbox{ for some } a \in
\Omega)\end{align*}
\noindent   where $\Omega$ is the set of those vectors $a =
(a_1,\ldots,a_n) \in \Z^n \backslash \{0\}$ which have at most
$n/2$ non-zero entries. But in \cite[Section 3.1]{kks} it is shown
that
$$ \P( a \cdot X_1 = \ldots = a \cdot X_n = 0 \hbox{ for some } a \in \Omega) \leq (1 + o(1))n^2 2^{-n}$$
and the claim follows.  (In fact, in the estimate in \cite{kks} one can permit as many as $n - 3 \log_2 n$ of the co-efficients
$a_1,\ldots,a_n$ to be non-zero).
\end{proof}

We will supplement these two easy lemmas with the following, somewhat more difficult, result.

\begin{proposition}[Medium combinatorial dimension estimate]\label{medium-dim}  Let $0 < \ep_0 \ll 1$,
and let $d_\pm \in D$ be such that $(\frac{3}{4} + 2\ep_0)^n < 2^{d_\pm-n} < \frac{100}{\sqrt{n}}$. Then we have
$$\sum_{V \in \Gr(d_\pm)} \P( A_V ) \leq o(1)^n,$$
where the rate of decay in the $o(1)$ quantity depends on $\ep_0$ (but not on $d_\pm$).
\end{proposition}

Note that $D$ has cardinality $|D| = O(n^2)$.  Thus if we combine
this proposition with Lemma \ref{kks-small} (with $\alpha :=
\frac{3}{4} + 2\ep_0$) and Lemma \ref{kks-large}, we see that we can bound the left-hand side of \eqref{pn-spanned-2}
by
$$ n (\frac{3}{4} + 2\ep_0)^n + n^2 o(1)^n + (1 + o(1))n^2 2^{-n} = (\frac{3}{4} + 2\ep_0 + o(1))^n.$$
Since $\ep_0$ is arbitrary, Theorem \ref{main} follows.

It thus suffices to prove Proposition \ref{medium-dim}.  This we shall do in later sections, but for now let us just remark that this
argument suggests a way to obtain Conjecture \ref{conj1} in full.  A comparison of the three bounds in Lemma \ref{kks-small}, Lemma \ref{kks-large},
and Proposition \ref{medium-dim} reveals that it is the estimate in Lemma \ref{kks-small} which is by far the most inefficient.  Indeed, from the other two lemmas we now see that the conjecture $P_n = (\frac{1}{2} + o(1))^n$ is equivalent to proving that
$$ \sum_{d_\pm \in D: 2^{d_\pm-n} \leq (\frac{3}{4} + 2\ep_0)^n} \,\,\,\, \sum_{V \in \Gr(d_\pm)} \P( A_V )
\leq (\frac{1}{2} + o(1))^n$$
for at least one value of $\ep_0 > 0$.  This bound is sharp (except for the $o(1)$ factor) as can be seen by considering the event that two of
the $X_i$ are equal up to sign. On the other hand, if we let $B$ denote the event that $a_1 X_1 + \ldots + a_n X_n = 0$
for some $a_1,\ldots,a_n \in \Z^n \backslash \{0\}$ with at most $n - 3 \log_2 n$ non-zero entries, it is shown in \cite[Section 3.1]{kks} that
$\P(B) = (1 + o(1))n^2 2^{-n}$ (this is basically the transpose of Lemma \ref{kks-large}).
Thus it would suffice to show
$$ \sum_{d_\pm \in D: 2^{d_\pm-n} \leq (\frac{3}{4} + 2\ep_0)^n}
\,\,\,\,  \sum_{V \in \Gr(d_\pm)} \P( A_V \backslash B) =
(\frac{1}{2} + o(1))^n.$$ Based on Proposition \ref{medium-dim}, we may
tentatively conjecture that in fact we have the stronger statement
\begin{equation}\label{conjecture-strong}
 \sum_{d_\pm \in D: 2^{d_\pm-n} \leq (\frac{3}{4} + 2\ep_0)^n}
 \,\,\,\, \sum_{V \in \Gr(d_\pm)} \P( A_V \backslash B)
= o(1)^n;
\end{equation}
this for instance would even imply the stronger conjecture in Conjecture \ref{conj1}.  But our Fourier-based methods seem to hit a natural limit at
$2^{d_\pm-n} \sim (3/4)^n$ and so are unable to obtain any further progress towards \eqref{conjecture-strong}.

\section{The general approach}\label{general-sec}

We now informally discuss the proof of Proposition \ref{medium-dim}; the rigorous proof will begin in Section \ref{excep-reduction}.
We start with the trivial bound
\begin{equation}\label{trivial-bound}
 \sum_{V \in \Gr(d_\pm)} \P( A_V ) \le 1
\end{equation}
that arises simply because any vectors $X_1,\ldots,X_n$ can span at most one space $V$. To improve upon this trivial bound,
the key innovation in \cite{kks} is to replace $X$ by another
random variable $Y$ which tends to be more concentrated on
subspaces $V$ than $X$ is.  Roughly speaking, one seeks the
property
\begin{equation}\label{c-compare}
 \P( X \in V ) \leq c \P( Y \in V)
\end{equation}
for some absolute constant $0 < c < 1$ and for all (or almost all)
subspaces $V \in \Gr(d_\pm)$.  From this property, one expects (heuristically, at least)
\begin{equation}\label{c-compare-iter}
 \P(A_V) = \P( X_1,\ldots,X_n \hbox{ span } V ) \leq c^n \P( Y_1,\ldots,Y_n \hbox{ span } V ),
\end{equation}
where $Y_1,\ldots,Y_n$ are iid samples of $Y$, and then by
applying the trivial bound \eqref{trivial-bound} with $Y$ instead
of $X$, we would then obtain a bound of the form $\sum_{V \in \Gr(d_\pm)} \P( A_V ) \leq c^n$,
at least in principle.  Clearly, it will be desirable to make $c$ as small as possible; if we can make $c$
arbitrarily small, we will have established Proposition \ref{medium-dim}.

The random variable $Y$ can be described as follows.  Let $0 \leq
\mu \leq 1$ be a small absolute constant (in \cite{kks} the value
$\mu = \frac{1}{108} e^{-1/108}$ was chosen), and let
$\eta^{(\mu)}$ be a random variable taking values in $\{-1,0,1\}
\subset F$ which equals 0 with probability $1-\mu$ and equals $+1$
or $-1$ with probability $\mu/2$ each.  Then let $Y :=
(\eta^{(\mu)}_1,\ldots,\eta^{(\mu)}_n) \in F^n$, where
$\eta^{(\mu)}_1,\ldots,\eta^{(\mu)}_n$ are iid samples of
$\eta^{(\mu)}$.  By using some Fourier-analytic arguments of
Hal\'asz, a bound of the form
$$ \P( X \in V) \leq C \sqrt{\mu} \P(Y \in V)$$
was shown in \cite{kks}, where $C$ was an absolute constant (independent of
$\mu$), and $V$ was a hyperplane which was \emph{non-degenerate}
in the sense that its combinatorial dimension was not too close to $n$.  For $\mu$ sufficiently small,
one then obtains \eqref{c-compare} for some $0 < c < 1$, although one cannot make $c$ arbitrarily small without
shrinking $\mu$ also.

There are however some technical difficulties with this approach,
arising when one tries to pass from \eqref{c-compare} to
\eqref{c-compare-iter}.  The first problem is that the random
variable $Y$, when conditioned on the event $Y \in V$, may
concentrate on a lower dimensional subspace on $V$, making it
unlikely that $Y_1,\ldots,Y_n$ will span $V$.  In particular, $Y$
has a probability of $(1-\mu)^n$ of being the zero vector, which
basically means that one cannot hope to exploit \eqref{c-compare}
in any non-trivial way once $\P(X \in V) \leq (1-\mu)^n$. However,
in this case $V$ has very low combinatorial dimension and
Lemma \ref{kks-small} already gives an exponential gain.

Even when $(1-\mu)^n < \P(X \in V) \leq 1$, it turns out that it is still
not particularly easy to obtain \eqref{c-compare-iter}, but one
can obtain an acceptable substitute for this estimate by only
replacing some of the $X_j$ by $Y_j$. Specifically,  one can try
to obtain an estimate roughly of the form
\begin{equation}\label{gammab} \P( X_1,\ldots,X_n \hbox{ span } V ) \leq c^{m}
\P( Y_1,\ldots,Y_m, X_1, \ldots, X_{n-m} \hbox{
span } V )
\end{equation}
where $m$ is equal to a suitably small multiple of $n$ (we will eventually take $m \approx n/100$).  Strictly speaking,
 we will also have to absorb an additional ``entropy'' loss of $\binom{n}{m}$ for technical reasons, though as we will be taking $c$
 arbitrarily small, this loss will ultimately be irrelevant.

The above approach (with some minor modifications) was carried out
rigorously in \cite{kks} to give the bound $P_n = O(.999^n)$ which
has been  improved slightly to $O(.939^n)$ in \cite{TV1}. There
are two main reasons why the final gain in the base was relatively
small. Firstly, the chosen value of $\mu$ was small (so the
$n(1-\mu)^n$ error was sizeable), and secondly the value of $c$
obtained was relatively large (so the gain of $c^n$ or
$c^{(1-\gamma)n}$ was relatively weak).  Unfortunately, increasing
$\mu$ also causes $c$ to increase, and so even after optimizing
$\mu$ and $c$ one falls well short of the conjectured bound.

In this paper we partially resolve this problem as follows.  To reduce all the other losses to $(\frac{3}{4} + 2\ep_0)^n$ for some small $\eps_0$,
we increase $\mu$ up to $1/4 - \ep_0/100$, at which point the arguments of Hal\'asz
and \cite{kks}, \cite{TV1} give \eqref{c-compare} with $c = 1$. The value
$1/4$ for $\mu$ is optimal as it is the largest number satisfying
the pointwise inequality
$$ |\cos(x)| \leq  (1-\mu) + \mu \cos(2x) \hbox{ for all } x \in \R,$$
which is the Fourier-analytic analogue of \eqref{c-compare} (with $c=1$).
At first glance, the fact that $c=1$ seems to remove any utility to \eqref{c-compare}, as the above argument relied
on obtaining gains of the form $c^n$ or $c^{(1-\gamma)n}$. However, we can proceed further by subdividing the
collection of hyperplanes $\Gr(d_\pm)$ into two classes, namely the
\emph{unexceptional} spaces $V$ for which
$$ \P( X \in V) < \eps_1 \P(Y \in V)$$
for some small constant $0 < \eps_1 \ll 1$ to be chosen later (it will be much smaller than $\eps_0$), and the
\emph{exceptional} spaces for which
\begin{equation}\label{V-exceptional}
\eps_1 \P(Y \in V) \leq \P(X \in V) \leq \P(Y \in V).
\end{equation}
The contribution of the unexceptional spaces can be dealt with by
the preceding arguments to obtain a very small contribution (at
most $\delta^n$ for any fixed $\delta>0$ given that we set $\eps_1
=\eps_1(\gamma,\delta)$ suitably small), so it remains to consider
the exceptional spaces $V$.  The key new observation (which uses
the Fourier analytic methods of Hal\'asz and \cite{kks}) is that
the condition \eqref{V-exceptional} can be viewed as a statement
that the ``spectrum'' of $V$ (which we will define later) has
``small doubling constant''.  Once one sees this, one can apply
theorems from inverse additive number theory (in particular a
variant Freiman's theorem) to obtain some strong structural
control on $V$ (indeed, we obtain what is essentially a complete
description of the exceptional spaces $V$, up to constants
depending on $\eps_0$ and $\eps_1$).  This then allows us to
obtain a fairly accurate bound for the total number of exceptional
spaces $V$, which can then be used to estimate their contribution
to \eqref{pn-spanned} satisfactorily.

\section{Reduction to exceptional spaces}\label{excep-reduction}

We now rigorously carry out the strategy outlined in Section \ref{general-sec}.
We first pick the parameter $\mu$ as
\begin{equation}\label{mu-def}
\mu := \frac{1}{4} - \frac{\eps_0}{100},
\end{equation}
and let $Y \in \{-1,0,1\}^n \subset F^n$ be the random variables defined in Section \ref{general-sec}
using using this value of $\mu$.  We let $Y_1,\ldots,Y_n$ be iid
samples of $Y$, independent of the previous samples $X_1,\ldots,X_n$ of $X$.
Thus we may write
\begin{equation}\label{XY-def}
X = \sum_{j=1}^n \eta^{(1)}_j e_j; \quad Y = \sum_{j=1}^n
\eta^{(\mu)}_j e_j
\end{equation}
where the $\eta^{(\mu)}_j$ are iid samples of the random variable
$\eta^{(\mu)}$ introduced in the previous section, $\eta^{(1)}_j$
are iid random signs, and $e_1,\ldots,e_n$ is the standard basis of $F^n$.

Next, we introduce a small parameter $\ep_1 > 0$ (which will be much smaller than $\ep_0$).
Let us call a space $V \in \Gr(d_\pm)$ \emph{exceptional} if
\begin{equation}\label{exceptional-def}
\P( X \in V) \geq \eps_1 \P( Y \in V )
\end{equation}
and \emph{unexceptional} if
\begin{equation}\label{unexceptional-def}
\P( X \in V) < \eps_1 \P( Y \in V ).
\end{equation}
As the terminology suggests, the exceptional spaces will turn out to be relatively rare compared to the unexceptional spaces, as
will be seen by comparing Lemma \ref{unex-estimate} and Lemma \ref{ex-estimate} below.

Proposition \ref{medium-dim} is now a consequence of the following two sub-lemmas.

\begin{lemma}[Unexceptional space estimate]\label{unex-estimate} We have
\begin{equation}\label{unexceptional-est}
\sum_{V \in \Gr(d_\pm): V \hbox{ unexceptional}} \P(A_V)
\leq
 2^{o(n)} 2^n \eps_1^{\eps_0 n/100}.
\end{equation}
where the decay rate in the $o(n)$ term can depend on $\ep_0, \ep_1$.
\end{lemma}

\begin{lemma}[Exceptional space estimate]\label{ex-estimate} We have
$$
\sum_{V \in \Gr(d_\pm): V \hbox{ exceptional}} \P(A_V)
\leq n^{- \frac{n}{2} + o(n)}
$$
where the decay rate in the $o(n)$ term can depend on $\ep_0, \ep_1$.
\end{lemma}

Indeed, upon combining Lemma \ref{unex-estimate} and Lemma \ref{ex-estimate} we obtain
$$
\sum_{V \in \Gr(d_\pm)} \P(A_V)
\leq 2^{o(n)} 2^n \eps_1^{\eps_0 n/100} + n^{-n/2 + o(n)},
$$
which implies Proposition \ref{medium-dim}
since $\eps_1$ can be chosen arbitrarily small.

It remains to prove Lemma \ref{unex-estimate} and Lemma \ref{ex-estimate}.  We prove Lemma \ref{unex-estimate}, which is simpler
in this section, and leave the more difficult Lemma \ref{ex-estimate} to later sections.

\begin{proof}[Proof of Lemma \ref{unex-estimate}]
We adapt some arguments form \cite{kks}, \cite{TV1}; the estimates here are somewhat cruder (and simpler) than in those papers, because we have the additional factor of $\eps_1$ which can absorb several losses which can arise in this argument.  We begin with a lemma.

\begin{lemma}[Weighted Odlyzko lemma]\label{Y-odd}\cite{kks}  Let $0 \leq d \leq n$. If $W$ is any $d$-dimensional subspace of $F^n$, then
$$ \P( Y \in W ) \leq (1-\mu)^{n - d}.$$
\end{lemma}

\begin{proof}  Write $Y = (\eta_1,\ldots,\eta_n)$.  Then
there exists some $d$-tuple $(j_1,\ldots,j_d)$ of
co-ordinates which determine the element of $W$; this means that
once $\eta_{j_1},\ldots,\eta_{j_d}$ are selected, there is only
one possible choice for each of the remaining $n-d$ independent random variables $\eta_i$ if one
wishes $Y$ to lie in $W$.  Since any given value is attained by an
$\eta_i$ with probability at most $1-\mu$, the claim follows.
\end{proof}

Let $m$ be the nearest integer to $\eps_0 n/100$.
In addition to the random variables $X_1,\ldots,X_n$, we
create more random variables $Y_1,\ldots,Y_{m}$ which have the same distribution as $Y$, and which are independent of each other and of $X_1,\ldots,X_n$.

Using Lemma \ref{Y-odd}, we have

\begin{lemma} Let $V \in \Gr(d_\pm)$ be unexceptional, and let $B_{V,m}$ denote the event
that $Y_1,\ldots,Y_m$ are linearly independent and lie in $V$.  Then
$$ \P( B_{V,m} ) \geq 2^{-o(n)} (2^{d_\pm - n} / \eps_1)^m.$$
\end{lemma}

\begin{proof} Using Bayes' identity, we can factorize
$$ \P(B_{V,m}) = \prod_{i=1}^m \P(B_{V,i}|B_{V,i-1})$$
where $B_{V,i}$ is the event that $Y_1,\ldots,Y_i$ are linearly independent and lie in $V$.  Let $W_i$ be the linear subspace of $F^n$ spanned by
$Y_1,\ldots,Y_i$.  Then conditioning on any fixed value of $Y_1,\ldots,Y_{i-1}$ in the event $B_{V,i-1}$, we have
$$ \P(B_{V,i}|B_{V,i-1}) = \P(Y \in V) - \P(Y \in W_i),$$
since $Y_i$ is independent of $Y_1,\ldots,Y_{i-1}$ and has the same distribution as $i$.  By \eqref{unexceptional-def}, \eqref{discrete-dim} we have
$$ \P(Y \in V) \geq \frac{1}{\eps_1} \P(X \in V) \geq \frac{1}{\eps_1} 2^{-1/n} 2^{d_\pm - n}.$$
On the other hand, since $W_i$ has dimension at most $m$, we see from Lemma \ref{Y-odd}, \eqref{mu-def} that
$$ \P(Y \in W_i) \leq (1 - \mu)^{n-m} = (\frac{3}{4} + \frac{\eps_0}{100})^{n-m}.$$
Since $m = \eps_0 n/100 + O(1)$, and $(\frac{3}{4} + 2\ep_0)^n < 2^{d_\pm-n}$ by hypothesis, we thus see that
$$ \P(Y \in W_i) = O(1/n) \P(Y \in V)$$
(for instance) if $n$ is sufficiently large depending on $\eps_0$.  Thus we have
$$ \P(B_{V,i}|B_{V,i-1}) \geq (1 - O(1/n)) 2^{d_\pm - n} / \eps_1.$$
Multiplying this together for all $1 \leq i \leq m$, the claim follows.
\end{proof}

To apply this lemma, observe that $A_V$ and $B_{V,m}$ are clearly independent, as they involve independent sets of random variables.
Thus we have
$$ \P(A_V) \leq 2^{o(n)} (2^{d_\pm - n} / \eps_1)^{-m} \P( A_V \wedge B_V ).$$
Now let $X_1,\ldots,X_n,Y_1,\ldots,Y_m$ be in the event $A_V \wedge B_V$.  Since $Y_1,\ldots,Y_m$ are linearly independent in $V$, and
$X_1,\ldots,X_n$ span $V$, we see that there exist $n-m$ vectors in $X_1,\ldots,X_n$ which, together with $Y_1,\ldots,Y_m$, span $V$.
The number of possibilities for such vectors is $\binom{n}{n-m}$, which we can crudely bound by $2^n$.  By symmetry we thus have
$$ \P(A_V) \leq 2^{o(n)} (2^{d_\pm - n} / \eps_1)^{-m} 2^n \P( C_V ),$$
where $C_V$ is the event that $Y_1,\ldots,Y_m, X_1,\ldots,X_{n-m}$ span $V$, and that $X_{n-m+1},\ldots,X_n$ lie in $V$.  By independence
we then have
$$ \P(A_V) \leq 2^{o(n)} (2^{d_\pm - n} / \eps_1)^{-m} 2^n \P( Y_1,\ldots,Y_m,X_1,\ldots,X_{n-m} \hbox{ span} V ) \P(X \in V)^{m}.$$
Applying \eqref{discrete-dim} we conclude
$$ \P(A_V) \leq 2^{o(n)} (2^{d_\pm - n} / \eps_1)^{-m} 2^n \P( Y_1,\ldots,Y_m,X_1,\ldots,X_{n-m} \hbox{ span} V ) (2^{d_\pm -n})^{m}.$$
Summing this over $V$, and observing that each collection of vectors $Y_1,\ldots,Y_m,X_1,\ldots,X_{n-m}$ can only span a single space $V$,
we conclude
$$
\sum_{V \in \Gr(d_\pm): V \hbox{ unexceptional}} \P(A_V)
\leq 2^{o(n)} (2^{d_\pm - n} / \eps_1)^{-m} 2^n (2^{d_\pm -n})^{m}
$$
and the claim follows from the choice of $m$.
\end{proof}

\section{Structure of exceptional subspaces}

To conclude the proof of Theorem \ref{main}, the only remaining
task is to prove Lemma \ref{ex-estimate}.  Observe from \eqref{discrete-dim} that
$$ \P(A_V) \leq \P(X \in V)^n \leq 2^{-n(n-d_\pm)}$$
for all $V \in \Gr(d_\pm)$.  Thus it will suffice to show that
\begin{equation}\label{exceptional-est}
 |\{ V \in \Gr(d_\pm): V \hbox{ exceptional} \}| \leq n^{- \frac{n}{2} + o(n)}
 2^{n(n-d_\pm)}.
\end{equation}
To do this we require a certain structural theorem
about exceptional subspaces, which is the main novelty in our
argument, and the one which requires tools from inverse additive
number theory such as Freiman's theorem.

We need some notation before we can describe the structure theorem.

\begin{definition}[Progressions]
Let $M_1, \dots, M_r$ be integers and let $a, v_1, \dots, v_r$ be
non-zero elements of $F$.  The set
$$ P := \{ a + m_1 v_1 + \ldots + m_r v_r: -M_j/2 < m_j < M_j/2 \hbox{ for all } 1 \leq j \leq r \}$$
is called a \emph{generalized arithmetic progression} of rank $r$;
we say that the progression is \emph{symmetric} if $a=0$. We say
that $P$ is \emph{proper} if  the map $$(m_1,\ldots,m_r) \mapsto
m_1 v_1 + \ldots + m_r v_r$$ is injective when $-M_j/2 < m_j <
M_j/2$. If $P$ is proper and symmetric, we define the \emph{$P$-norm} $\|v\|_P$
of a point $v = m_1 v_1 + \ldots + m_r v_r$ in $P$ by the formula
$$ \| m_1 v_1 + \ldots + m_r v_r \|_P := \left(\sum_{i=1}^r (\frac{|m_i|}{M_i})^2\right)^{1/2}.$$
\end{definition}

Let $V \in \Gr(d_\pm)$ be an exceptional space, with a
representation of the form
\begin{equation}\label{V-rep}
 V = \{ (x_1,\ldots,x_n) \in F^n: x_1 a_1 + \ldots + x_n a_n = 0 \}
 \end{equation}
for some elements $a_1, \ldots, a_n \in F$. We shall refer to
$a_1,\ldots,a_n$ as the \emph{defining co-ordinates} for $V$.

\begin{theorem}[Structure theorem]\label{structure}  There is a constant $C=C(\ep_0, \ep_1)$
 such that the following holds. Let $V$ be an exceptional hyperplane in $
 \Gr(d_\pm)$ and $a_1,\ldots,a_n$ be its defining co-ordinates.
Then  there exist integers

\begin{equation}\label{rank-bound}
1 \leq r \leq C
\end{equation}
 and $M_1,\ldots,M_r \geq 1$ with the volume bound
\begin{equation}\label{m-max}
 M_1 \ldots M_r \leq C 2^{n - d_\pm}
\end{equation}
and non-zero elements $v_1, \ldots, v_r \in F$ such that the
following holds.

\begin{itemize}

\item (i) (Defining coordinates lie in a progression) The
symmetric generalized arithmetic progression

$$ P := \{ m_1 v_1 + \ldots + m_r v_r: -M_j/2 < m_j < M_j/2 \hbox{ for all } 1 \leq j \leq r \}$$
is proper and contains all the $a_i$.

\item (ii) (Bounded norm)  The $a_i$ have small $P$-norm:
\begin{equation}\label{p-mush}
 \sum_{j=1}^n \| a_j \|_P^2 \leq C
\end{equation}

\item (iii) (Rational commensurability) The  set $\{
v_1,\ldots,v_r \} \cup \{a_1,\ldots,a_n\}$ is contained in the set
\begin{equation}\label{vr-rank}
\{ \frac{p}{q} v_1 : p, q \in \Z; q \neq 0; |p|, |q| \leq n^{o(n)}
\}.
\end{equation}

\end{itemize}

\end{theorem}

\begin{remark} The condition (i) asserts that the defining co-ordinates
$a_1,\ldots,a_n$ of $V$ lie in a fairly small generalized arithmetic
progression with bounded rank $O(1)$. The condition \eqref{vr-rank}
asserts, furthermore, that this progression is contained in a (rank one) arithmetic
progression of length $n^{o(n)}$.  Thus we have placed $\{a_1,\ldots,a_n\}$
inside two progressions, one of small size but moderately large rank, and one of rank one but of fairly huge size (but the dimensions
are still smaller than $n^{n/2}$ or $|F|$).
\end{remark}

\begin{remark} The structure theorem is fairly efficient; if
$a_1,\ldots,a_n$ obey the conclusions of the theorem, then from
\eqref{p-mush} and the theory of random walks (taking advantage of
the boundedness of $r$, which will play the role of dimension) we
expect the random variables $\eta^{(1)}_1 a_1 + \ldots +
\eta^{(1)}_1 a_n$ and $\eta^{(\mu)}_1 a_1 + \ldots +
\eta^{(\mu)}_1 a_n$ to be distributed fairly uniformly on $P$ and
decay rapidly away from $P$, and so we expect the probabilities
$\P( X \in V)$ and $\P(Y \in V)$ to be comparable (up to a
constant).
\end{remark}

We shall prove the structure theorem in later sections.
 For the
remainder of this section, we show how the structure theorem can
be used to prove Lemma \ref{ex-estimate}.  In the sequel all of the $o()$ factors
are allowed to depend on $\ep_0$, $\ep_1$.

Let $V \in \Gr(d_\pm)$ be an exceptional hyperplane.
Then $V$ has $|F|-1$ representations of the form \eqref{V-rep}, one for each non-zero normal vector of $V$.
Let us call an $n$-tuple $(a_1,\ldots,a_n)$ of elements in $F$
\emph{exceptional} if it obeys the conclusions of the structure
theorem for at least one progression $P$, thus we have
\begin{align*}
  |\{ V \in \Gr(d_\pm): V \hbox{ exceptional} \}|
= &\frac{1}{|F|-1} \\
|\{ (a_1,\ldots,a_n) \in &F^n: (a_1,\ldots,a_n) \hbox{
exceptional} \}|,
\end{align*}
and so it now suffices to show

\begin{align*} &\,\,\,\,  |\{ (a_1,\ldots,a_n) \in F^n: (a_1,\ldots,a_n)
\hbox{ exceptional} \}|  \\ &\leq n^{o(n)} n^{-n/2}
2^{n(n-d_\pm)} |F|. \end{align*}

We now need to count the number of exceptional $(a_1,\ldots,a_n)$.
We first observe that we may fix the parameter $r$ in the
conclusion of the structure theorem, since the number of $r$ is at
most  $C = n^{o(n)}$.  Similarly we may fix each $M_1,\ldots,M_r$,
since the total number of choices here is at most $(C 2^{n -
d_\pm})^r = n^{o(n)}$.

The number of possible choices for $v_1$ is at most $|F|$.  Once this vector is selected,
we see from \eqref{vr-rank} that there are at most $n^{o(n)}$
possible choices for each of the remaining
$r-1 = O(1)$ vectors.  Putting the estimates together, we can conclude that
the total number of possible vectors $(v_1,\ldots,v_r)$ which
could be chosen is at most $n^{o(n)} |F|$.

It now suffices to show that, for
each fixed choice of $r$, $M_1,\ldots,M_r$, $v_1,\ldots,v_r$
(which in particular fixes $P$), that the number of possible
exceptional $(a_1,\ldots,a_n)$ is at most $n^{o(n)} n^{-n/2}
2^{n(n-d_\pm)}$. We are going to use the bound on the $P$-norms
of the $a_i$. Our task is to show that for any constant $C$

\begin{equation}\label{nanda}
\Big|\{ (a_1,\ldots,a_n) \in P^n: \sum_{j=1}^n \| a_j \|_P^2 \leq
C \} \Big| \leq n^{o(n)} n^{-n/2} 2^{n(n-d_\pm)} .
\end{equation}
We shall use Gaussian-type methods.  To start, notice that

$$ \Big| \{ (a_1,\ldots,a_n) \in P^n: \sum_{j=1}^n \| a_j \|_P^2 \le C  \} \Big|  \\
\leq n^{o(n)} \sum_{(a_1,\ldots,a_n) \in P^n} \exp(- n
\sum_{j=1}^n \| a_j \|_P^2). $$

\noindent The right hand side can be rewritten as

$$n^{o(n)} \Big( \sum_{a \in P} \exp(- n \| a \|_P^2) \Big)^n =
n^{o(n)} \Big(\prod_{j=1}^r \sum_{|m_j| \leq M_j} e^{-n m_j^2 /
M_j^2} \Big)^n. $$

\noindent Using the elementary bound
$$ \sum_{|m| \leq M} e^{-n m^2 /
M^2} = O(1+ M n^{-1/2} )$$

\noindent we can bound the right hand side from
above by

$$  n^{o(n)} \big( \prod_{1 \leq j \leq r}  (1 + \frac{M_j}{\sqrt{n}}) \big)^n
. $$

\noindent Since $r = O(1)$ and $M_1,\ldots,M_r \geq 1$, the product $\prod_{1 \leq j \leq r}  (1 +
\frac{M_j}{\sqrt{n}})$ can be bounded (rather crudely) by $O( 1 + n^{-1/2} M_1 \ldots
M_r)$. This leads us to
 the following bound

$$ n^{o(n)}  (1 + n^{-1/2} M_1 \ldots M_r)^n. $$

\noindent Now we use the information (see (\ref{m-max}) in
Structure Theorem) that $M_1 \ldots M_r$ is at most $C
2^{n-d_\pm}$. The bound becomes
$$ n^{o(n)} ( 1 + n^{-1/2} C 2^{n-d_\pm} )^n =  n^{-n/2+o(n)}
2^{n(n-d_\pm)}$$
as desired, thanks to the hypothesis $2^{d_\pm-n} \leq 100 / \sqrt{n}$.

This concludes the proof of Theorem \ref{main} except for the
Structure Theorem, which we turn to next. This proof requires a
heavy use of tools from additive combinatorics. We are going to
introduce these tools in the next section.

\section{Generalized arithmetic progressions and inverse
theorems}

In this section we review some notation and results from additive combinatorics, and briefly sketch how
they will be used to prove Theorem \ref{structure}.  The rigorous proof of this theorem will begin in Section \ref{halasz}.

Let $A$  and $B$ be  finite subsets of $G$. By $A+B := \{ a+b: a \in A, b \in B\}$ we
denote the
set of all elements of $G$ which can be represented as the sum of
an element from $A$ and an element from $B$; $A-B := \{a-b: a \in A, b\in B\}$ is defined
similarly. Moreover, for any positive integer $l$, we recursively define
$lA:=(l-1)A+A$, thus $lA$ is the set of $l$-fold sums in $A$ (allowing repetition).

The \emph{doubling constant} of a finite non-empty set $A \subseteq F$ is
defined to be the the ratio $|A+A|/|A|$. It is easy
to see that if $A$ is a dense subset (with constant
density $\delta$) of a proper generalized arithmetic progression of constant
rank $r$, then the doubling number of $A$ is bounded from above by
a constant depending on $\delta$ and $r$ (in fact it is bounded
by $\delta^{-1} 2^r$).

In the mid 1970s, Freiman \cite{Fre} proved a remarkable result in the converse
direction,
which asserts that being a dense subset  of a proper arithmetic
progression of constant rank $r$ is the only reason for a finite
set $A$ of integers to have small doubling constant:

\begin{theorem}[Freiman's theorem]\label{Freiman}\cite{Fre} For any constant $C$ there are
constant $r$ and $\delta$ such that the following holds. For any
finite set $A$ of integers such that $|A+A| \le C |A|$, there is a
proper arithmetic progression $P$ of rank $r$ such that $A \subseteq
P$ and $|A|/|P| \ge \delta$. \end{theorem}

Freiman's proof was rather complex. A cleaner proof, based on
Freiman's ideas, is presented by Bilu in \cite{Bil}. Ruzsa
\cite{Ruz} gave a different and quite short proof; this was then refined by
Chang \cite{chang}.  There are explicit bounds known on the values of $r$ and $\delta$ in terms of
$C$, but for our application these bounds will only influence the $o(1)$ term in our final result and so
we will not keep track of them here.

Freiman's theorem was initially phrased in the integers $\Z$, but it can easily be transferred to a finite field $F$ of prime order
if $F$ is sufficiently large:

\begin{theorem}\label{GR}  For any constant $C$ there are
constant $r$ and $\delta$ such that the following holds. Let $F$ be a finite field of prime order, and let $A$ be a non-empty subset of $F$ such
that  $|A+A|\leq C|A|$. Then, if $|F|$ is sufficiently large depending on $|A|$, there is a generalized arithmetic progression $P$ of
rank $r$ such that $A\subset P$ and $|A|/|P| \ge \delta$.
\end{theorem}

 This theorem follows from Freiman's original theorem and  ``rectification'' theorems such as \cite[Theorem 3.1]{lev} (see
also \cite{green-ruzsa3}) to map $A$ via a \emph{Freiman
isomorphism} to the the integers $\Z$ (using the hypothesis that
$F$ is large depending on $A$; $|F| \geq \exp(C |A|)$ would do).
It is also a special case of the version of Freiman's theorem
established in \cite{green-ruzsa4} for an arbitrary abelian group.
In fact, that result allows one to remove the hypothesis that
$|F|$ is sufficiently large depending on $|A|$.

A generalized arithmetic progression is not always proper. We can,
however, make this assumption whenever we like (at the cost of a
constant factor) thanks to the following lemma.

\begin{lemma}[Progressions lie inside proper progressions]\label{proper} There is a constant $C$ such that the following holds.
Let $P$ be a generalized arithmetic progression of rank $r$
in an abelian group $G$.  Suppose that every non-zero element of $G$ has order at least $r^{Cr^3} |P|$.
Then there exists a proper arithmetic progression $Q$ of rank at most $r$ containing $P$ and
$$ |Q| \leq r^{Cr^3} |P|.$$
\end{lemma}

The proof of this lemma arises purely from the geometry of numbers
(in particular, Minkowski's second theorem) and is independent from the rest of this paper. In order
not to distract the reader, we defer this rather technical proof in the
Appendix.

Finally, we are going to need the following lemma, based on a
covering  argument  of Ruzsa \cite{ruzsa-group}. We say that a set
$A$ is \emph{symmetric} if the set $-A :=\{-a, a \in A\}$ is equal to $A$.

\begin{lemma}[Sum set estimates]\label{Ruzsa} Let $A$ be a symmetric finite subset of an
abelian group $G$ such that $|4A| \le C|A|$ for some $C \geq 1$. Then for any $k \ge
4$
$$|kA| \le \binom{ {C +k-3} }{ {k-2}} C|A|. $$
\end{lemma}

\begin{proof}
 We can assume that $k \ge 4$.  Consider
the sets $\eta + A$ as $\eta$ ranges inside $3A$. Each set has
cardinality $|A|$, and is contained inside $4A$. Thus we may find
a maximal disjoint collection $\{ \xi + A: \eta \in X \}$ where $X
\subseteq 3A$ has cardinality $|X| \leq |4A|/|A| \leq C$. Then for
any $\xi \in 3A$ there exists $\eta \in X$ such that $\eta +
\Lambda$ intersects $\xi + A$, otherwise this would contradict
maximality. But this implies that $\xi \in X + A - A = 2A + X$.
Thus we have
$$ A + 2A = 3A \subseteq 2A + X.$$

Iterating this we obtain
$$ kA \subset 2A + (k-2)X$$
for all $k \geq 2$.  Thus

$$ |kA| \leq |2A| |(k-2)X| \leq C|A|(k-2)X|. $$

To conclude the proof, notice that for any $l$

$$|lX| \le {\binom{ |X| +l-1 } {l}   }. $$

\end{proof}

In fact one can replace the condition $|4A| \le C|A|$ by  a weaker
condition $|3A| \le C|A|$.

We are going to use Theorem \ref{GR} to derive a structural
property of the defining vectors $\{a_1, \dots, a_n \}$ of $V$
in  Structure Theorem. To be more precise, we are going to apply
Theorem \ref{GR} to a set $A$ with small doubling constant
which contains the majority of the defining vectors $\{a_1, \dots, a_n \}$.
Thus we obtain a generalized arithmetic
progression  $P$ containing $A$ (and with it most of the $a_i$).
It will be easy to extend $P$ to contain all $a_i$ without violating
the properties we care about.

The application  of Theorem \ref{GR} is not immediate. In the next
section, we shall provide a Fourier-analytic argument which leads
to the definition of $A$. The starting point of this argument is
based on the method \cite{kks}, which is motivated by   an
observation of Hal\'asz \cite{Hal}; this method shall be the focus
of the next section.

\section{Hal\'asz-type arguments}\label{halasz}

Let $V \in \Gr(d_\pm)$ be an exceptional space with
the representation \eqref{V-rep}.  In order to obtain the desired
structural control on the defining co-ordinates $a_1,\ldots,a_n$, we shall first
use Fourier analysis to first gain some strong control on the
``spectrum'' $\Lambda$ of $a_1,\ldots,a_n$ (which we will define
in \eqref{lambda-def}), and then apply the inverse Fourier transform to
recover information about $a_1,\ldots,a_n$.

We turn to the details.
By Fourier expansion we have
$$ 1_{(x_1,\ldots,x_n) \in V} = \frac{1}{|F|} \sum_{\xi \in F} e_p( x_1 a_1 \xi + \ldots +  x_n a_n \xi )$$
for all $(x_1,\ldots,x_n) \in F^n$, where $p = |F|$ and $e_p$ is the primitive
character $e_p(x) := e^{2\pi i x/p}$.  By \eqref{XY-def} we thus
have (by linearity of expectation and by independence)
\begin{align*}
\P( X \in V) &= \E( 1_{X \in V} ) \\
&= \frac{1}{|F|}  \sum_{\xi \in F} \E( e_p( \eta^{(1)}_1 a_1 \xi + \ldots + \eta^{(1)}_n a_n \xi ) ) \\
&= \frac{1}{|F|}  \sum_{\xi \in F} \prod_{j=1}^n \E( e_p( \eta^{(1)}_j a_j \xi ) ) \\
&= \frac{1}{|F|} \sum_{\xi \in F} \prod_{j=1}^n \cos( 2 \pi a_j
\xi / p ) .
\end{align*}
In particular we have
\begin{align*}
\P( X \in V)
&\leq \frac{1}{|F|} \sum_{\xi \in F} \prod_{j=1}^n |\cos( 2 \pi a_j \xi / p )| \\
&= \frac{1}{|F|} \sum_{\xi \in F} \prod_{j=1}^n |\cos( \pi a_j \xi / p )|\\
&= \frac{1}{|F|} \sum_{\xi \in F} \prod_{j=1}^n (\frac{1}{2} +
\frac{1}{2} \cos( 2\pi a_j \xi / p ))^{1/2}
\end{align*}
where the first identity follows from the substitution $\xi
\mapsto \xi/2$ on $F$, noting that the quantity $|\cos(\pi a_j
\xi / p)|$ is still periodic in $\xi$ with period $p$.
Arguing similarly with $Y$, we have
\begin{align*}
\P( Y \in V) &= \frac{1}{|F|}  \sum_{\xi \in F} \prod_{j=1}^n \E( e_p( \eta^{(\mu)}_j a_j \xi ) ) \\
&= \frac{1}{|F|} \sum_{\xi \in F} \prod_{j=1}^n ((1-\mu) + \mu
\cos( 2 \pi a_j \xi / p )).
\end{align*}
To summarize, if we introduce the non-negative quantities
\begin{equation}\label{f-slip}
f(\xi) := \prod_{j=1}^n (\frac{1}{2} + \frac{1}{2} \cos( 2\pi a_j
 \xi / p ))^{1/2};  \quad g(\xi) := \prod_{j=1}^n ((1-\mu) +
\mu \cos( 2 \pi a_j \xi / p ))
\end{equation}
then we have
\begin{equation}\label{XY-form}
 \P(X \in V) \leq \frac{1}{|F|} \sum_{\xi \in F} f(\xi); \quad \P(Y \in V) = \frac{1}{|F|} \sum_{\xi \in F} g(\xi).
\end{equation}

We now give a crucial comparison estimate between $f$ and $g$.

\begin{lemma}\label{comparison} For all $\xi \in F$, we have $f(\xi) \leq g(\xi)^{1/4\mu}$.
\end{lemma}

\begin{proof}
By \eqref{f-slip} suffices to show that
$$ (\frac{1}{2} + \frac{1}{2} \cos \theta)^{1/2} \leq ((1-\mu) + \mu \cos\theta)^{1/4\mu}$$
for any $\theta$.  Writing $\cos \theta = 1 - 2x$ for some $0 < x
\leq 1$ (the $x=0$ case being trivial), this can be rearranged to
become
$$ \frac{\log(1 - x) - \log(1 - 0)}{x} \leq \frac{\log(1 - 2\mu x) - \log(1-0)}{2\mu x}.$$
But this follows from the concavity of the function $\log(1-x)$ on
$0 \leq x \leq 1$ and the fact (from \eqref{mu-def}) that $0 <
2\mu < 1$.
\end{proof}

Since $g(\xi)$ is clearly bounded by 1, and since $\mu < 1/4$ by
\eqref{mu-def}, it follows that for every $\xi$
\begin{equation}\label{fg}
f(\xi) \leq g(\xi),
\end{equation}
which when combined with \eqref{XY-form} gives
$$ \P( X \in V ) \leq \P( Y \in V ).$$

We now refine this argument to exploit the
additional (and critical)  hypothesis \eqref{exceptional-def} that
$$ \P( X \in V) \geq \eps_1 \P( Y \in V ). $$

Let $\ep_2 > 0$ be a sufficiently small positive constant (compared to
$\ep_1$). Define the \emph{spectrum} $\Lambda \subseteq F$ of
$\{a_1,\ldots,a_n\}$ to be the set
\begin{equation}\label{lambda-def}
 \Lambda := \{ \xi \in F: f(\xi) \geq \eps_2 \}
 \end{equation}
Note that $\Lambda$ is symmetric around the origin: $\Lambda =
-\Lambda$.  Next we make the elementary observation that
\begin{equation}\label{elementary}
1 - 100 \|x\|^2 \leq \cos(2\pi x) \leq 1 - \frac{1}{100} \|x\|^2
\end{equation}
(say), where $\|x\|$ is the distance of $x$ to the nearest integer. From \eqref{f-slip} we thus have
(with room to spare)
$$ f(\xi) \leq \exp( - \frac{1}{1000}  \sum_{j=1}^n \| a_j \cdot \xi / p \|^2 )$$
and hence there is a constant $C(\ep_2)$ depending on $\ep_2$ such
that

\begin{equation}\label{lambda-ineq} (\sum_{j=1}^n \| a_j \cdot \xi / p \|^2 )^{1/2} \leq C(\eps_2)
\end{equation}
for all $\xi \in \Lambda$.

We now obtain some cardinality bounds on $\Lambda$ and of the iterated sumsets $k\Lambda$.

\begin{lemma}\label{lambda-bounds}  There is a constant $C$
depending on $\ep_0, \ep_1, \ep_2$ such that
\begin{equation}\label{lambda-b}
C^{-1} 2^{-(n-d_\pm)} |F| \leq |\Lambda| \leq C 2^{-(n-d_\pm)} |F|
\end{equation}
 Furthermore, for every integer $k \ge 4$
\begin{equation}\label{k-bound}
 |k\Lambda| \leq { \binom{C +k-3}{k-2}} C 2^{-(n-d_\pm)}
 |F|.
 \end{equation}

\end{lemma}

\begin{remark} The lemma implies that $\Lambda$ has a small
doubling constant. At this point one could apply Theorem \ref{GR} to
gain further control on $\Lambda$. While this can be done, it is
more convenient for us to apply the inverse Fourier transform and
work on a certain ``dual'' set to $\Lambda$, which then contains
most of the $a_j$.  A key point in \eqref{k-bound}, which we will
exploit in the proof of Lemma \ref{A-double} below, is that the
growth of constants is only polynomial in $k$ rather than
exponential.

Since $\ep_i$ depends on $\ep_0, \dots, \ep_{i-1}$, it would
suffice to say that $C$ depends only on $\ep_0$.
\end{remark}

\begin{proof}  From \eqref{exceptional-def}, \eqref{XY-form} we have
$$ \frac{1}{|F|} \sum_{\xi \in F} f(\xi) \geq \eps_1 \frac{1}{|F|} \sum_{\xi \in F} g(\xi).$$
But from Lemma \ref{comparison}, the definition of $\Lambda$, and
the crucial fact that $\mu < 1/4$, we know that
\begin{align*}
\frac{1}{|F|} \sum_{\xi \not \in \Lambda} f(\xi)
&\leq \eps_2^{1 - 4\mu}  \frac{1}{|F|} \sum_{\xi \not \in \Lambda} f^{4\mu}(\xi)\\
&\leq \eps_2^{1-4\mu} \frac{1}{|F|} \sum_{\xi \in F} g(\xi).
\end{align*}

Now we are going to make an essential use of  the assumption that
$\mu$ is less than $1/4$. Given this assumption,  we can
guarantee that the contribution of $\xi$ outside the spectrum is
negligible, by choosing $\eps_2$ sufficiently compared to
$\eps_1$. This results in the following equalities

$$\sum_{\xi \in \Lambda} f(\xi) =\Theta (\sum_{\xi \in F} f(\xi))
=\Theta (\sum_{\xi \in F} g(\xi)) =\Theta (|F| \P (X \in V)) =
\Theta (2^{d_{\pm} -n} |F|), $$

\noindent where the constants in the $\Theta()$ notation depend on $\ep_0,
\ep_1, \ep_2$. The bounds on $|\Lambda|$ now follows directly from
the fact that for any $\xi \in \Lambda$, $\ep_2 \le f(\xi) \le1$.

It remains to prove \eqref{k-bound}.   In view of Lemma
\ref{Ruzsa}, it suffices to show that there is a constant $C$ such
that

$$|4\Lambda | \le C |\Lambda|. $$

Notice that if $\xi \in 4\Lambda$, then by \eqref{lambda-ineq} and
the triangle inequality we have
$$ (\sum_{j=1}^n \| a_j \cdot \xi / p \|^2 )^{1/2} \leq C(\eps_2),$$

\noindent for some constant $C(\ep_2)$ depending on $\ep_2$ (we
abuse the notation a little bit, as the two $C(\ep_2)$ are not
necessarily the same). From \eqref{elementary}, \eqref{f-slip} we conclude a bound of the form

$$ f(\xi) \geq c(\eps_2) > 0.$$

Thus

$$ |4\Lambda| \leq c(\eps_2)^{-1}  \sum_{\xi \in G} f(\xi) \le  c(\eps_2)^{-1}  \sum_{\xi \in F} g(\xi) = \Theta (2^{d_{\pm} -n} |F|)
=\Theta (|\Lambda|) $$

\noindent concluding the proof.

\end{proof}

We now pass from control of the spectrum back to control on
$\{a_1,\ldots,a_n\}$, using the inverse Fourier transform. For any
$x \in F$, define the norm $\|x\|_\Lambda$ by
$$ \|x\|_\Lambda := \Big( \frac{1}{|\Lambda|^2} \sum_{\xi,\xi' \in \Lambda} \| x (\xi-\xi') / p \|^2 \Big)^{1/2}.$$
It is easy to see that this is a quantity between 0 and 1 which
obeys the triangle inequality $\|x+y\|_\Lambda \leq \|x\|_\Lambda
+ \|y\|_\Lambda$. From the triangle inequality again we have
\begin{align*}
 \|x\|_\Lambda &\leq ( \frac{1}{|\Lambda|^2} \sum_{\xi,\xi' \in \Lambda} \| x  \xi / p \|^2)^{1/2}
+ ( \frac{1}{|\Lambda|^2} \sum_{\xi,\xi' \in \Lambda} \| x \xi' / p \|^2)^{1/2}\\
&= 2 ( \frac{1}{|\Lambda|} \sum_{\xi} \| x  \xi / p
\|^2)^{1/2}.
\end{align*}
Thus by square-summing \eqref{lambda-ineq} for all $\xi \in
\Lambda$, we obtain
\begin{equation}\label{A-manage}
 \sum_{j=1}^n \| a_j \|_\Lambda^2 \leq C(\eps_2).
\end{equation}
Thus we expect many of the $a_j$ to have small $\Lambda$ norm.  On
the other hand, we are going to show that the set of elements with
small $\Lambda$ norm has constant doubling number, thanks to the
following lemma.

\begin{lemma}\label{A-double}  There is a constant $C$ such that the following holds.
Let $A \subseteq F$ denote the ``Bohr set''
$$ A := |\{ x \in F: \|x\|_\Lambda \leq \frac{1}{100} \}|.$$
Then we have
$$ C^{-1} 2^{n-d_\pm} \leq |A| \leq |A+A| \leq C2^{n-d_\pm}.$$
\end{lemma}

\begin{proof}  Let us first establish the upper bound.  We introduce the normalized Fourier transform of $\Lambda$:
$$ h(x) := \frac{1}{|\Lambda|} \sum_{\xi \in \Lambda} e_p( x \xi ).$$
If $x \in A+A$, then we have $\|x\|_\Lambda \leq \frac{1}{50}$ by
the triangle inequality. 
Using  \eqref{elementary}, we conclude


\begin{align*}
|h(x)|^2 &= \Re h(x) \overline{h(x)} \\
&= \Re \frac{1}{|\Lambda|^2} \sum_{\xi, \xi' \in \Lambda} e_p( x  (\xi-\xi') ) \\
&= \frac{1}{|\Lambda|^2} \sum_{\xi, \xi' \in \Lambda} \cos( 2 \pi x (\xi-\xi') / p ) \\
&\ge \frac{1}{|\Lambda|^2} \sum_{\xi, \xi' \in \Lambda} 1 - 100 \| x  (\xi-\xi') / p \|^2 \\
&= 1 - 100 \| x \|_\Lambda^2.
\end{align*}

\noindent As  $\|x\|_\Lambda \leq \frac{1}{50}$, it follows that

$$|h(x) | \ge \sqrt {1 - 100 (1/50)^2 } > 1/2. $$

On the other hand, from the Parseval identity we have
\begin{equation}\label{pars}
 \sum_{x \in F} |h(x)|^2 = |F|/|\Lambda|;
\end{equation}
combining these two estimates we obtain
$$ |A+A| \leq 4 |F|/|\Lambda|$$
and the claim now follows from \eqref{lambda-b}.

Since the bound $|A| \leq |A+A|$ is trivial, it now suffices to
prove the lower bound on $A$. This proof has two parts. First we
show that if $|h(x)|$ is large (close to 1), then $\| x
\|_\Lambda$ is small (close to zero). Next, we prove that there
are many $x$ such that  $|h(x)|$ is large.

The first part is simple and similar to the above argument. Using
the fact that  $\|x\|_\Lambda \leq \frac{1}{50}$ and
\eqref{elementary}, we have

\begin{align*}
|h(x)|^2 &= \Re h(x) \overline{h(x)} \\
&= \Re \frac{1}{|\Lambda|^2} \sum_{\xi, \xi' \in \Lambda} e_p( x  (\xi-\xi') ) \\
&= \frac{1}{|\Lambda|^2} \sum_{\xi, \xi' \in \Lambda} \cos( 2 \pi x (\xi-\xi') / p ) \\
&\leq \frac{1}{|\Lambda|^2} \sum_{\xi, \xi' \in \Lambda} 1 - \frac{1}{100} \| x  (\xi-\xi') / p \|^2 \\
&= 1 - \frac{1}{100} \| x \|_\Lambda^2.
\end{align*}
It follows that if $|h(x)| \ge 1-10^{-4}$ then $x \in A$.

Now we
are going to show that the set of those $x$ such that $|h(x)| \ge
1-10^{-4}$ has cardinality $\Theta (2^{n-d_{\pm}})$.
 This proof is
trickier and we need to  adapt an argument from
\cite{green-ruzsa3} (see also \cite{green-ruzsa4}), which requires the full strength of
\eqref{k-bound}. Let $k$ be a large integer to be chosen later,
and for each $\xi \in F$, let $r_k(\xi)$ be the number of
representations of the form $\xi = \xi_1 + \ldots + \xi_k$ where
$\xi_1, \ldots, \xi_k \in \Lambda$. Clearly we have
$$\sum_{\xi \in k\Lambda} r_k(\xi) = |\Lambda|^k$$
and so by Cauchy-Schwarz
$$ \sum_{\xi \in F} |r_k(\xi)|^2 \geq |\Lambda|^{2k} / |k\Lambda|.$$
By Plancherel we have
$$ \sum_{\xi \in F} |r_k(\xi)|^2 = \frac{1}{|F|} \sum_{x \in F} |\sum_{\xi \in F} r_k(\xi) e_p(x \xi)|^2.$$
But we have
\begin{align*}
|\sum_{\xi \in F} r_k(\xi) e_p(x  \xi)| &=
|\sum_{\xi_1,\ldots,\xi_k \in \Lambda} e_p(x (\xi_1+\ldots+\xi_k)| \\
&= |\Lambda|^k |h(x)|^k
\end{align*}
and hence
$$ \sum_{x \in F} |h(x)|^{2k} \geq |F| / |k\Lambda|.$$
Combining this with \eqref{pars} we obtain
$$ \sum_{x \in F: |h(x)| \geq (|\Lambda|/2|k\Lambda|)^{1/(2k-2)} } |h(x)|^{2k}
\geq |F|/2|k\Lambda|;$$ since $h(x)$ is trivially bounded by 1, we
thus conclude
$$ |\{ x \in F: |h(x)| \geq (|\Lambda|/2|k\Lambda|)^{1/(2k-2)} \}| \geq |F|/2|k\Lambda|.$$
Inserting the bounds \eqref{lambda-b}, \eqref{k-bound} we have

$$ (|\Lambda|/2|k\Lambda|)^{1/(2k-2)} \ge [{ \binom{C +k-3}
{k-2}}C]^{1/(2k-2)} $$

\noindent for some constant $C$. Choose $k$ sufficiently large
(say $k=C^2 +10^{10} $) we can guarantee that the right hand side
is at least $1-10^{-4}$. On the other hand $k$ is still a
constant, so $|k\Lambda| =O( 2^{d_{\pm} -n} |F|)$. It follows that

$$ |\{ x \in F: |h(x)| \geq 1 - 10^{-4} \}| =\Theta( 2^{n-d_\pm}),$$

\noindent proving the lower bound on $A$.
\end{proof}

\section {Proof of Structure Theorem}

We are now in position to prove the Structure Theorem.  In this section we allow the constants
in the $O()$ notation to depend on $\ep_0, \ep_1, \ep_2$.
From Lemma \ref{A-double} we
already know that $A$ has doubling constant $O(1)$. Furthermore, it is
easy to show that $A$ contains most of the $a_i$. Indeed, if $a_i
$ is not in $A$ then its $\Lambda$-norm is at least $1/100$;
\eqref{A-manage} then shows that
$$|\{ 1 \leq i \leq n: a_i \not \in A \}| = O(1).$$

We now apply Theorem \ref{GR} to $A$, to conclude that there is a  generalized
arithmetic progression  $P$ of rank $O(1)$ containing $A$ and
with $|A|= \Theta (|P|) = \Theta (2^{n-d_{\pm}})$. By Lemma
\ref{proper} we can assume that $P$ is proper (note that every
non-zero element of $F$ has order $|F| = p \ge n^{n/2}$, which is
certainly much larger than $|A|$).

This progression $P$ is close to, but not quite, what we need for
the structure theorem.  Thus we shall now perform a number of
minor alterations to $P$ in order to make it exactly match the
conclusions of the structure theorem.  Each of our alterations
will preserve the following properties claimed in the structure theorem, namely
\begin{equation} \label{rankP} \hbox{rank}(P) = O(1) \hbox{ and } \end{equation}
\begin{equation} \label{volumeP} |P| =M_1 \dots M_r = O(2^{n-d_{\pm}}), \end{equation}
though of course the $O()$ constants may worsen with each of the alterations.

First, we can include the elements of $\{a_1, \dots, a_n\}
\backslash A$ in $P$ simply by adding each exceptional $a_j$ as a
new basis vector $v'_j$ (with the corresponding length $M'_j$ set
equal to 3). This enlargement of $P$ (which by abuse of notation
we shall continue to call $P$) now contains $A$ and all of the
$a_j$, and obeys the bounds
\eqref{rankP}, \eqref{volumeP} with slightly worse constants. Notice
that by  Lemma \ref{proper} we can always assume $P$ is proper, without affecting the validity
of \eqref{rankP} and \eqref{volumeP} other than in the constants.

It now remains to verify the estimate \eqref{p-mush} on the
$P$-norm, as well as the rational commensurability claim.

To verify \eqref{p-mush}, consider the sum $\sum_j \|a_j\|_P^2 $.
For the exceptional $a_j$ (which were added to $P$ afterwards)
$\|a_j\|_P=O(1)$. However, the number of these $a_j$ is $O(1)$. For
the remaining $a_j$, notice that $P$  contains $ka_j$ for all  $k
< \frac{1}{100 \|a_j\|_{\Lambda}}$. Thus, \eqref{p-mush} follows
from \eqref{A-manage} which asserts that $\sum_j
\|a_j\|_{\Lambda}^2 =O(1)$. This leaves the verification of the rational commensurability
as the only remaining task.

Let us, as usual, write
$$P= \{ m_1v_1+ \dots +m_r v_r| -M_j/2 \le m_j \le M_j/2 \}. $$

\noindent We consider $P$ together with the map $\Phi: P
\rightarrow \R^r$ which maps $m_1v_1+ \dots +m_r v_r$ to $(m_1,
\dots, m_r)$. Since $P$ is proper, this map is bijective. We set $U
:=\{a_1, \dots, a_n\} \subseteq P$.

 We know that the progression $P$ contains
$U$, but we do not know yet that $U$ ``spans'' $P$ in the sense
that the set $\Phi(U)$ has full rank (i.e., it spans $\R^r$).
However, this is easily rectified by a rank reduction argument.
Suppose that $\Phi(U)$ did not have full rank, we are going to
produce a new proper generalized arithmetic progression $P'$ which
still contains $U$ and satisfies \eqref{volumeP} but has rank
strictly smaller than the rank of $P$. Since the rank of $P$ is $O(1)$, the process must terminate after
at most $O(1)$ iterations.

We are going to produce $P'$ as follows. If $\Phi(U)$ does not
have full rank, then it is  contained in a hyperplane of $\R^r$.
In other words, there exist integers $\alpha_1,\ldots,\alpha_r$
whose greatest common divisor is one and   $\alpha_1 m_1 + \ldots
+ \alpha_r m_r = 0$ for all $(m_1,\ldots,m_r) \in \Phi(U)$. Let
$w$ be an arbitrary element of $F$ and replace the each basis
vector $v_j$ in $P$ by $v_j - \alpha_j w$. The new progression $P'$ will
continue to contain $U$, since we have
$$ m_1 (v_1 - \alpha_1 w) + \ldots + m_r (v_r - \alpha_r v_r) = m_1 v_1 + \ldots + m_r v_r $$

\noindent for all $(m_1, \dots, m_r) \in \Phi(U) $. We can assume,
without loss of generality, that $\alpha_r $  is not divisible by
$p$, the characteristic of $F$ (at least one of the $\alpha_j$
must be so). Select $w$ so that

$$v_r-\alpha_r w  =0. $$

This shows that $P'$ has rank $r-1$. Moreover, its volume is $M_1,
\dots M_{r-1}$ which is at most the volume of $P$. We can use
Lemma \ref{proper} to guarantee that $P'$ is proper without
increasing the rank. Furthermore, it is easy to see that
\eqref{p-mush} is still valid with respect to $P'$. The argument
is completed. Thus, from now on we can assume that $\Phi(U)$ spans
$\R^r$. This information will be useful later on.

 Define a \emph{highly rational} number to be a
number (in $\Q$ or $F$) which is of the form $a/b$ where $a,b$ are
integers such that $a,b = n^{o(n)}$ and $b \neq 0$. Define a
\emph{highly rational linear combination} of vectors
$w_1,\ldots,w_k$ to be any expression of the form $q_1 w_1 +
\ldots + q_k w_k$ where $q_1,\ldots,q_k$ are highly rational.

We say that a set $W$ (of vectors)  {\it economically spans} a set
$U$ (of vectors) if there are numbers $a,b= n^{o(n)}$ such that
each $u \in U$ can be represented as a highly rational linear
combination of vectors in $W$ where the numerators and
denominators of the coefficient are uniformly  bounded from above
by $a$ and $b$, respectively.

The moral of our arguments below is the following: with respect to a
bounded number of operations, any linear algebraic statement which
holds for ``span'' also holds for ``economically span''. This relies
on the fact that  any algebraic expression involving at most
$O(1)$ highly rational numbers will again produce a highly
rational number. In other words, the highly rational numbers behave heuristically like
a subfield of $F$.

Let us now turn to the details.  First we  make some simple
remarks about the notion of economically spanning. It is a
transitive in the sense that if $W$ economically spans $ U$ and
$U$ economically span $ T$, then $W$ economically spans $T$ (but
with slightly worse $o(n)$ constants). Furthermore, it is closed
under union, in the sense that if $W$ economically spans $ X$ and
$Y$,  then it economically spans the union of $X$ and $Y$. Of
course we are going to use these rules only a bounded number
$O(1)$ of times.

Consider the set of vectors $U := \{a_1,\ldots,a_n\}$.  We know
that there is a set of vectors of cardinality $r$ which
economically span $U$, namely $\{v_1,\ldots,v_r\}$. Now we are
going to make use of the fact that $\Phi(U)$ spans $\R^r$. Since
each vector in $\Phi(U)$ has coordinates at most $\max_{1 \leq i
\leq n} M_i \le 2^n = n^{o(n)}$ and $r=O(1)$, it follows from
Cramer's rule that $\Phi(U)$ economically spans the basic vectors
$\{e_1, \dots, e_r\}$. As the map $\Phi$ is bijective ($P$ is
proper), this means $U$ economically spans $\{v_1, \dots, v_r\}$.

We now claim
that in fact there is a  single element $v_1$ which economically
spans $U$. This clearly implies rational commensurability.

 Let $s$ be the smallest number such that there is a
subset of size $s$ of $v_1,...v_r$ economically spanning $U$.
Without loss of generality, we assume that $v_1,...,v_s$ spans $U$
and write

\begin{equation} \label{equa:ai} a_i = \sum_{j=1}^s c_{ij} v_j, \end{equation} where $c_{ij}$ are highly
rational. Let us now consider two cases:

\begin{itemize}
\item The matrix $C=(c_{ij})$ has rank one. In this case the numbers $a_i/a_j$
are highly rational (as all the $c_{ij}$ are). But $\Phi(U)$ has
full rank in $ R^r$ and economically spans $v_1,..,v_r$, so
$v_i/v_j$ are also highly rational. This means that $v_1$ (or any
of the $v_i$, for that matter)  economically spans $v_1, \dots,
v_r$, and hence also $U$ by transitivity.

\item The matrix $C $ has rank larger than one. Recall that
$(a_1,..,a_n)$ is the normal vector of a hyperplane spanned by
$\pm 1$ vectors.  Thus,  there is a $\pm 1$ vector $w$ which is
orthogonal to $a= (a_1,...a_n)$, but not orthogonal to $C$.
Express the inner product of $a$ and $w$ as a linear combination
of $v_1,\dots, v_s$ using (\ref{equa:ai}). Observe that all of the
coefficients are highly rational. Moreover, the value of the
combination is 0, but there is at least one non-zero coefficient.
This implies that we can express one of $v_1, \dots, v_s $ as a
linear combination of the others with highly rational coefficients
and this contradicts the minimality of $s$.
\end{itemize}

The proof of the structure theorem is complete.
\endprf

\section {Appendix: Proof of Lemma \ref{proper} }

In this section we give a proof of Lemma \ref{proper}.  This result appears to be ``folklore'', being discovered independently by
Gowers and Walters (private communication) and Ruzsa (private communication), but does not appear to be explicitly in the
literature.

We first need some notation and tools from the geometry of numbers.
 The first tool is the following theorem of Mahler on bases of lattices, which is a variant
of Minkowski's second theorem (and is in fact proven using this theorem).

\begin{lemma}\label{wbasis}
Let $\Gamma$ be a lattice of full rank in $\R^d$ (i.e. a discrete additive subgroup of $\R^d$ with $d$ linearly independent generators).  Then there exists linearly independent vectors $w_1, \ldots, w_d \in \Gamma$
which generate $\Gamma$, and such that
\begin{equation}\label{w-magic}
|w_1| \ldots |w_d| \leq (Cd)^{Cd} {\rm mes}(\R^d/\Gamma)
\end{equation}
where ${\rm mes}(\R^d/\Gamma)$ is the volume of a fundamental domain of $\Gamma$ and $C > 0$ is an absolute constant.
\end{lemma}

\begin{proof}  (Sketch)  Note that the standard Minkowski basis $v_1,\ldots,v_d$ of $\Gamma$ (with respect
to the unit ball $B_d$) will obey \eqref{w-magic} but need not
generate $\Gamma$ (consider for instance the lattice generated by
$\Z^d$ and $(\frac{1}{2},\ldots,\frac{1}{2})$ for dimensions $d
\geq 5$).  However, a simple algorithm of Mahler gives a genuine
basis $w_1,\ldots,w_d$ of $\Gamma$ such that $w_j = t_{j,1} v_1 +
\ldots + t_{j,j} v_j$ for some real scalars $|t_{j,i}| \leq 1$,
which thus also obeys \eqref{w-magic} with some degradation in the
$(Cd)^{Cd}$ factor.  See, for instance \cite[Chapter 8]{cassels}
for more details.
\end{proof}

If $d \geq 1$ is an integer,
and $M = (M_1,\ldots,M_d)$ and $N = (N_1,\ldots,N_d)$ are elements of $\R^d$ such that $M_i \leq N_i$ for all $1 \leq i \leq d$,
we use $[M,N]$ to denote the discrete box
$$ [M,N] := \{ (n_1,\ldots,n_d) \in \Z^d: M_i \leq n_i \leq N_i \hbox{ for all } 1 \leq i \leq d \}.$$
Also, if $G$ is an additive group, $n = (n_1,\ldots,n_d) \in \Z^d$, and $v = (v_1,\ldots,v_d) \in G^d$, we use $n \cdot v \in G$
to denote the quantity $n \cdot v = n_1 v_1 + \ldots + n_d v_d$, where of course we use $nx := x + \ldots + x$ to denote the $n$-fold
sum of $x$ (with the obvious modification if $n$ is negative).  We also use $[M,N] \cdot v \subseteq G$ to denote the set
$[M,N] \cdot v := \{ n \cdot v: n \in [M,N]\}$.  Note that a progression of rank $d$
in $G$ is nothing more than a set of the form $a + [0,N] \cdot v$ for some $a \in G$, some $d$-tuple $N = (N_1,\ldots,N_d) \in \Z_+^d$,
some other $d$-tuple $v = (v_1,\ldots,v_d) \in G^d$.

We now give a ``discrete John's theorem'' which shows that the intersection of a convex symmetric body (i.e. a bounded open convex symmetric subset of $\R^d$) and a lattice of full rank is essentially equivalent to a progression.

\begin{lemma}[Discrete John's theorem]\label{djt}  Let $B$ be a convex symmetric body in $\R^d$, and let $\Gamma$ be a lattice in $\R^d$ of full rank.
Then there exists a $d$-tuple $$w = (w_1,\ldots,w_d) \in
\Gamma^d$$ of linearly independent vectors in $\Gamma$ and and a
$d$-tuple $N = (N_1,\ldots,N_d)$ of positive integers such that
$$ (-N,N) \cdot w \subseteq B \cap \Gamma \subseteq ( - d^{Cd} N, d^{Cd} N ) \cdot w$$
where $C$ is an absolute constant.
\end{lemma}

\begin{proof}
We first observe using John's theorem \cite{john} and an invertible linear transformation that we may assume
without loss of generality that $B_d \subseteq B \subseteq d \cdot B_d$, where $B_d$ is the unit ball in $\R^d$.
We may assume $d \geq 2$ since the claim is easy otherwise.

Now let $w = (w_1, \ldots, w_d)$ be as in Lemma \ref{wbasis}.  For each $j$, let $N_j$ be the least integer greater than $1/d|w_j|$.
Then from the triangle inequality we see that $|n_1 w_1 + \ldots + n_d w_d| < 1$ whenever $|n_j| < N_j$, and hence $(-N,N) \cdot w$ is contained
in $B_d$ and hence in $B$.

Now let $x \in B \cap \Gamma$.  Since $w$ generates $\Gamma$, we have $x = n_1 w_1 + \ldots + n_d w_d$ for some integers $n_1,\ldots,n_d$;
since $B \subseteq d \cdot B_d$, we have $|x| \leq d$.  Applying Cramer's rule to solve for $n_1,\ldots,n_d$ and \eqref{w-magic}, we have
\begin{align*}
 |n_j| &= \frac{|x \wedge w_1 \ldots w_{j-1} \wedge w_{j+1} \wedge w_d|}{|w_1 \wedge \ldots \wedge w_d|}\\
 &\leq \frac{|x| |w_1| \ldots |w_d|}{|w_j| |w_1 \wedge \ldots \wedge w_d|}\\
 &\leq d d^{Cd} / |w_j|\\
 &\leq d^{C' d} N_j,
\end{align*}
and hence $x \in ( - d^{Cd} N, d^{Cd} N ) \cdot w$ as desired.
\end{proof}

We can now prove Lemma \ref{proper}, which we restate here with slightly different notation.

\begin{lemma}\label{proper-torsion} There is a constant $C_0$ such that the following holds.
Let $P$ be a symmetric progression of rank $d$
in a abelian group $G$, such that every non-zero element of $G$ has order at least $d^{C_0 d^3} |P|$.  Then there exists a symmetric proper progression $Q$ of rank at most $d$ containing $P$ and
$$ |Q| \leq d^{C_0 d^3} |P|.$$
\end{lemma}

\begin{proof}
This claim is analogous to the basic linear algebra statement that
every linear space spanned by $d$ vectors is equal to a linear
space with a \emph{basis} of at most $d$ vectors.  Recall that the
proof of that linear algebra fact proceeds by a descent argument,
showing that if the $d$ spanning vectors were linearly dependent,
then one could exploit that dependence to ``drop rank'' and span
the same linear space with $d-1$ vectors.

We may assume that $P$ has the form $P = [-N/2,N/2] \cdot v$
for some $N = (N_1,\ldots,N_d)$ and $v = (v_1,\ldots,v_d)$.
We induct on $d$.  The case  $d=1$ is easy. Now suppose
inductively that $d\geq 2$, and the claim has already been proven
for $d-1$ (for arbitrary groups $G$ and arbitrary progressions
$P$).  Thanks to the induction hypothesis, we can assume that $P$
is non-proper and all $N_i$ are at least one. As $P$ is not proper, there is $n \neq n' \in
[-N/2, N/2]$ such that
 $$n \cdot v = n' \cdot v. $$
  Let $\Gamma_0 \subseteq
\Z^d$ denote the lattice $\{ m \in \Z^d: m \cdot v = 0\}$, then it
follows that  $\Gamma_0 \cap [-N,N]$ contains at least one
non-zero element, namely $n'-n$.

Let $m = (m_1, \ldots, m_d)$ be a non-zero irreducible element of
$\Gamma_0 \cap [-N,N]$ (i.e. $m/n \not \in \Gamma_0$ for any integer $n>1$). By definition
\begin{equation}\label{movie}
m \cdot v = m_1 \cdot v_1 + \ldots + m_d \cdot v_d = 0.
\end{equation}
We now claim that $m$ is irreducible in $\Z^d$ (i.e. that $m_1, \ldots, m_d$ have no
common divisor). Indeed, if $m$ factored as $m = n \tilde m$ for some $n > 1$ and $\tilde m \in \Z^d$,
then $\tilde m \cdot v$ would be a non-zero element of $G$ of order $n \leq |P|$,
contradicting the hypothesis.

Our plan is to  contain $P$ inside a symmetric progression $Q$ of rank $d-1$
and cardinality
\begin{equation}\label{bomu}
|Q| \leq d^{Cd^2} |P|.
\end{equation}

If we can achieve this, then by the induction hypothesis we can
contain $Q$ inside a proper symmetric progression $R$ of rank at most $d-1$
and cardinality
$$ |R| \leq (d-1)^{C_0(d-1)^3} (Cd)^{Cd^2} |P|.$$

If $C_0$ is sufficiently large (but still independent of $d$) , then
the right-hand side is at most $d^{C_0d^3} |P|$, and we have closed
the induction hypothesis.

It remains to cover $P$ by a symmetric progression of rank at most
$d-1$ with the bound \eqref{bomu}.  Observe that $m$ lies in
$[-N,N]$, so the rational numbers $$m_1/N_1, \ldots, m_d/N_d$$ lie
between $-1$ and $1$.  Without loss of generality we may assume
that $m_d/N_d$ has the largest magnitude, namely
\begin{equation}\label{mp-max}
|m_d|/N_d \geq |m_i|/N_i
\end{equation}
for all $1 \leq i \leq d$.  By replacing $v_d$ with $-v_d$ if
necessary, we may also assume that $m_d$ is positive.

To exploit the cancellation in \eqref{movie} we introduce the
rational vector $q \in \frac{1}{m_d} \cdot \Z^{d-1}$ by the
formula
$$ q := (-\frac{m_1}{m_d}, \ldots, -\frac{m_{d-1}}{m_d}).$$
Since $\gcd(m_1, \dots, m_d)=1$, we see that for any integer $n$,
that $n \cdot q$ lies in $\Z^{d-1}$ if and only if $n$ is a
multiple of $m_d$.

Next, let $\Gamma \subset \R^{d-1}$ denote the lattice $\Gamma :=
\Z^{d-1} + \Z \cdot q$. Since $q$ is rational, this is indeed a
full rank lattice. We define the homomorphism $f: \Gamma
\rightarrow  Z$ by the formula
$$ f( (n_1,\ldots,n_{d-1}) + n_d q ) := (n_1,\ldots, n_d) \cdot v.$$
the condition \eqref{movie} ensures that this homomorphism is
indeed well defined, in the sense that different representations
$w = (n_1,\ldots,n_{d-1}) + n_d q$ of the same vector $w \in
\Gamma$ give the same value of $f(w)$. We also let $B \subseteq
\R^{d-1}$ denote the convex symmetric body
$$ B := \{ (t_1,\ldots,t_{d-1}) \in \R^{d-1}: -3N_j < t_j < 3N_j \hbox{ for all } 1 \leq j \leq d-1 \}.$$
We now claim the inclusions
$$ P \subseteq f( B \cap \Gamma ) \subseteq 5P - 5P.$$
To see the first inclusion, let $n \cdot v \in P$ for some $n \in
[0,N]$, then we have $n \cdot v = f( (n_1,\ldots,n_{d-1}) + n_d q
)$. From \eqref{mp-max} we see that the $j^{th}$ co-efficient of
$(n_1,\ldots,n_{d-1}) + n_d q$ has magnitude at most $3N_j$, and
thus $n \cdot v$ lies in $f(B \cap \Gamma)$ as claimed.  To see
the second inclusion, let $(n_1,\ldots,n_{d-1}) + n_d q$ be an
element of $B \cap \Gamma$.  By subtracting an integer multiple of
$m_d$ from $n_d$ if necessary (and thus adding integer multiples
of $m_1,\ldots,m_{d-1}$ to $n_1,\ldots,n_{d-1}$) we may assume
that $|n_d| \leq |m_d|/2$. By \eqref{mp-max} and the definition of
$B$, this forces $|n_j| \leq 5 N_j$ for all $1 \leq j \leq d$, and
hence
$$f((n_1,\ldots,n_{d-1}) + n_d q) = (n_1,\ldots,n_d) \cdot v \subseteq [-5N,5N] \cdot v = 5 P - 5P.$$
Next, we apply Theorem \ref{djt} to find vectors $w_1, \ldots,
w_{d-1} \in \Gamma$ and $M_1,\ldots,M_{d-1}$ such that
$$ (-M,M) \cdot w \subseteq B \cap \Gamma \subseteq ( - d^{Cd} M, d^{Cd} M ) \cdot w.$$
Applying the homomorphism $f$, we obtain
$$ (-M,M) \cdot f(w) \subseteq f(B \cap \Gamma) \subseteq ( - d^{Cd} M, d^{Cd} M ) \cdot f(w)$$
where $f(w) := (f(w_1),\ldots,f(w_{d-1})$.  Observe that $( -
d^{Cd} M, d^{Cd} M ) \cdot f(w)$ is a symmetric progression of rank $d-1$
which contains $f(B \cap \Gamma)$ and hence contains $P$.
Also, since $( - d^{Cd} M, d^{Cd} M ) \cdot f(w)$ can be covered by $O( (Cd)^{Cd^2} )$ translates
of $(-M, M) \cdot f(w)$, and similarly $5P - 5P$ can be covered by $O(C^d)$ translates of $P$, we have
\begin{align*}
|( - d^{Cd} M, d^{Cd} M ) \cdot f(w)|
&\leq (Cd)^{Cd^2} |f(B \cap \Gamma)|\\
&\leq(Cd)^{Cd^2} |5P-5P|\\
&\leq (Cd)^{Cd^2)} C^d |P|
\end{align*}
which proves \eqref{bomu}.  This completes the induction and
proves the theorem.
\end{proof}

\begin{remark} It is not hard to remove the hypothesis that $G$ is ``nearly torsion free'' in the sense that the order
of every non-zero element is much larger than $P$, however one must then drop the conclusion that $Q$ is symmetric (this can be seen simply by considering the case $G = \Z/2\Z$).
\end{remark}

\end{document}